\documentclass[11pt]{article}
\usepackage{amsmath,amssymb,amstext,amsthm}
\usepackage[dvips]{graphics}
 \numberwithin{equation}{section}
\numberwithin{thm}{section}

\numberwithin{lma}{section}
\newcommand{\ov}{\overline}
\usepackage{amsmath,amssymb,amstext,amsthm}
\usepackage[dvips]{graphics}
\usepackage{graphicx}
 \numberwithin{equation}{section}
\numberwithin{thm}{section}

\numberwithin{lma}{section} \numberwithin{Corollary}{section}

\begin{document}
\title{Asymptotic generalized bivariate extreme with random index}
\author{{\large M. A. Abd Elgawad}$^{a,b}$, A. M. Elsawah$^{a,c,d,}$\footnote{Corresponding author. E-mail: amelsawah@uic.edu.hk,
a.elsawah@zu.edu.eg, a\_elsawah85@yahoo.com}$~$, Hong Qin$^{a }$  and  Ting Yan$^{a}$ \\
{\footnotesize $^a$ {\it Faculty of Mathematics and Statistics,
Central China Normal University, Wuhan 430079, China}}\\
{\footnotesize $^b$  {\it Department of Mathematics, Faculty of
Science, Benha  University, Benha 13518, Egypt}}\\
{\footnotesize $^c$  {\it Department of Mathematics, Faculty of
Science, Zagazig University, Zagazig 44519, Egypt}}\\{\footnotesize $^d$ {\it Division of Science and Technology, BNU-HKBU United International College, Zhuhai 519085, China}}
}
\date{}
\maketitle
\begin{abstract}
In many biological, agricultural, military
activity problems and in some quality control problems, it is almost impossible to
have a fixed sample size, because some observations are always lost for various reasons.
Therefore, the sample size itself is considered frequently to be a random variable
(rv). The class of limit distribution functions (df's) of the random bivariate extreme generalized order statistics (GOS) from independent and identically distributed rv's are  fully characterized. When the random sample size is assumed to be independent of the basic variables and its df is assumed to converge weakly to a non-degenerate limit, the necessary and sufficient
conditions for the weak convergence of the random bivariate extreme GOS  are obtained. Furthermore, when the interrelation of the random size and the basic rv's is not
restricted, sufficient conditions of the convergence and the forms of the limit df's are deduced. Illustrative examples are given which lend further support to our theoretical results.
\end{abstract}
\section{ Introduction}
The concept of generalized order statistics (GOS) have been introduced by
Kamps (1995). It's enable a unified approach to  ascendingly
ordered random variables (rv's) as ordinary
order statistics (oos), sequential order statistics (sos), order statistics with non integral sample
size, progressively type II censored order statistics (pos), record values, $k$th record
values and Pfeifer's records. Let $\gamma_n=k>0,$ $\gamma_r=k+n-r+\sum_{j=r}^{n-1}m_j>0,$ $r=1,2,...,n-1,$ and $\tilde{m}=(m_1,\,m_2,\cdots\,,m_{n-1})\in \Re^{n-1}.$ Then the rv's
$X_{r:n}\equiv X(r,n,\tilde{m},k),r=1,2,...,n,$ are called GOS based on the
distribution function (df) $F$ with density function $f$ which are defined by their probability density function (pdf) \begin{equation*}
\begin{split}
 f_{1,2,...,n:n}^{(\tilde{m},k)}(x_1,x_2,...,x_n)=&\left(\prod_{j=1}^{n}\gamma_j\right)
\left(\prod_{j=1}^{n-1}(1-F(x_j))^{\gamma_j-\gamma_{j+1}-1}f(x_j)\right)\\
&\times(1-F(x_n))^{\gamma_n-1}f(x_n),
\end{split}
\end{equation*}
where $F^{-1}(0)\leq x_1\leq...\leq x_n\leq F^{-1}(1).$

In this work, we consider a wide subclass of GOS , by assuming $\gamma_{j}-\gamma_{j+1}=m+1>0.$  This subclass is known as $m-$GOS.
 Clearly  many important practical models of $m-$GOS are included such as oos, order statistics with non integer sample size and
sos. The marginal df's of the $r$th and $\grave{r}$th $m-$GOS (Nasri- Roudsari, 1996 and  Barakat, 2007)  are represented by  $\Phi^{(m,k)}_{r:n}(x)=I_{L_m(x)}(r,N-r+1)$ and $\Phi^{({m} ,k)} _{\grave{r}:n}(x)=
I_{L_m(x)}(N-R_r+1,R_r),$ respectively, where $\grave{r}=n-r+1, L_m(x) = 1 - (1 - F(x))^{m+1}, I_x(a, b)=\frac{1}{\beta(a,
b)}\int_{0}^{x} t^{a-1}(1 - t)^{b-1 }dt$  denotes the
incomplete beta ratio function,  $ N=\ell+n-1, R_r=\ell+r-1$ and $\ell=\frac{k}{m+1}.$   Moreover, by
using the results of Kamps (1995), we can write
explicitly the joint df's of the $r$th and $s$th $m-$GOS,  $m\neq-1,1\leq r<s\leq n,$  as:
$$\Phi_{r,s:n}^{(m,k)}( x,y)=
C_n^\star\int_{0}^{F(x)}\int_{\xi}^{F(y)}{\ov{\xi}^{m}\ov{\eta}^{\gamma_{s}-1}}(1-\ov{\xi}^{m+1})^{r-1}
(\ov{\xi}^{m+1}-\ov{\eta}^{m+1})^{s-r-1} d\eta d\xi,~x\leq y,\eqno (1.1)$$
where $C_n^\star=\frac{(m+1)^2\Gamma(N+1)}
{\Gamma(N-s+1)(r-1)!(s-r-1)!}$   and $\Gamma(.)$ is the usual gamma function. Recently Barakat et al. (2014a) studied the
limit df's of joint extreme $m-$GOS, for a fixed sample size. Moreover, the asymptotic behavior for bivariate df of the
lower-lower (l-l), upper-upper (u-u) and lower-upper (l-u) extreme $m-$GOS in Barakat et al. (2014b).

In the last few years much efforts had been devoted to investigate the limit df's of independent
rv's with random sample size. The appearance of this trend is naturally because
many applications require the consideration of such problem. For example, in many biological, agricultural and in some quality control problems, it is almost
impossible to have a fixed sample size because some observations always get lost
for various reasons. Therefore, the sample size $n$ itself is considered frequently to be a rv $\nu_n,$  where $\nu_n$ is independent of the basic variables (i.e., the original
random sample) or in some applications the interrelation of the basic variables  and the random sample size is not restricted.
Limit theorems for extremes with random sample size indexes have been thoroughly
studied in the above mentioned two particular cases :
 \begin{enumerate}
      \item The  basic variables and sample size index are independents
(see, Barakat, 1997).
      \item The interrelation of the basic variables and the random sample size is not restricted (see,
 Barakat and El Shandidy, 1990, Barakat, 1997 and Barakat et al., 2015a).
 \end{enumerate}

Our aim in this paper is to characterize the asymptotic behavior of the bivariate df's of the (l-l), (u-u) and (l-u) extreme $m-$ GOS  with random sample size.  When the random sample size is assumed to be independent of the basic variables and its df is assumed to converge weakly to a non-degenerate limit, the necessary and sufficient
conditions for the weak convergence of the random bivariate extreme $m-$GOS are obtained. Furthermore, when the interrelation of the random size and the basic rv's is not restricted, sufficient conditions of the convergence and the forms of the limit df's are deduced.   An illustrative examples are given which lend further support to our theoretical results. Throughout this paper the convergence in probability and the weak convergence, as $n\rightarrow\infty,$ respectively, denoted as
$"{\renewcommand{\arraystretch}{0.1}
\begin{array}{c}
\stackrel{\textstyle p}{\longrightarrow}\\
\scriptstyle{n}\end{array}}"$ and $"{\renewcommand{\arraystretch}{0.1}
\begin{array}{c}
\stackrel{\textstyle w}{\longrightarrow}\\
\scriptstyle{n}\end{array}}".$\\
\section{Asymptotic random bivariate extreme under $m-$GOS}
\subsection{Random sample size and basic rv's are independents }
In this subsection we deal with the weak convergence of bivariate  df's of the (u-u), (l-l) and (l-u) extreme $m-$GOS are fully characterized in Theorems 2.1, 2.2 and 2.3, respectively. When the sample size itself is a rv $\nu_n,$ which is assumed to be independent of the basic variables $X_{r:n},r=1,2,...,n.$\\
{\bf Theorem 2.1.} Consider  the following three conditions :
$$ \Phi_{\grave{r},\grave{s}:n}^{(m,k)}(x_{n}
,y_{n})=P({ Z}_{ \grave{r},\grave{s}:n}^{(n)}<\mathbf{x})=P(Z_{\grave{r}:n}^{ (n)}<x_{},Z_{{\grave{s}}:n}^{ (n)}<y){\renewcommand{\arraystretch}{0,1} \begin{array}{c}
\stackrel{w}{\longrightarrow}\\
\scriptstyle{n}\end{array}}\hat{\Phi}_{{r},{s}}^{(m,k)}(x_{}
,y_{}), x\leq y, \eqno (i)~$$
$$~ H_n(nz)= P(\frac{\nu_n}{n}<z){\renewcommand{\arraystretch}{0.1}
\begin{array}{c}
\stackrel{\textstyle w}{\longrightarrow}\\
\scriptstyle{n}\end{array}}H(z),~\eqno (ii)~$$
$$ {\Phi}_{\grave{r},\grave{s}:\nu_n}^{(m,k)}(x_{n}
,y_{n})=P( Z^{(n)}_{\grave{r},\grave{s}:\nu_{n} }<\mathbf{x})=P(Z_{\grave{r}:\nu_n}^{ (n)}<x_{},Z_{{\grave{s}}:\nu_n}^{ (n)}<y){\renewcommand{\arraystretch}{0,1} \begin{array}{c}
\stackrel{w}{\longrightarrow}\\
\scriptstyle{n}\end{array}}
    \hat{\Psi}_{{r},{s}}^{(m,k)}(x_{}
,y_{})$$
$$=\int_{0}^{\infty}  \overline{\Omega}_{r,s}^{(m,k)}(z\kappa_{1},z\kappa_{2})dH(z).~\eqno (iii)~$$
Then any two of the above conditions imply the remaining one, where $x_{n}=x_{1n}=a_n {x_{1}}+ b_n,y_{n}=x_{2n}=a_n {x_{2}}+ b_n, a_n>0, b_n $ are suitable normalizing constants, $\mathbf{x}=(x_{1},x_{2})=(x,y),\grave{r}=n-r+1<n-s+1=\grave{s},Z_{\grave{r_{i}}:n}^{(n)}=\frac{X^{}_{\grave{r_{i}}:n}-b_n}{a_n}, i=1,2,(\grave{r_{1}},\grave{r_{2}})=(\grave{r},\grave{s}),\hat{\Phi}_{{r},{s}}^{(m,k)}(x_{}
,y_{})$ is a non-degenerate df,  $~H(z)~$
is a df with $~H(+0)=0,$
$$ \overline{\Omega}_{r,s}^{(m,k)}(\kappa_{1},\kappa_{2})=\left\{\begin{array}{cc}
1-\Gamma_{R_{s}}(\kappa_{2}^{m+1}), ~~~~~~~~~~~~~~~~~~~~~~~~~~~~~~~~~~~~~~~~~~~~~~~~~~~~~~~~~~~~ ~ x\geq y,\\
1-\Gamma_{R_{r}}(\kappa_{1}^{m+1}) -\frac{1}{\Gamma(R_{r})}\int_{\kappa_{1}^{m+1}}^{\infty}I_{\left(\frac{\kappa_{2}^{m+1}}{u}\right)}(R_{s},R_{r}-R_{s})u^{R_{r}-1}{e^{-u}}du,  ~ x\leq y,\\
\end{array}
\right.$$
$\Gamma_{r}(x)=\frac{1}{\Gamma(r)}\int_{0}^{x}\theta^{r-1}e^{-\theta}d\theta$ denotes the incomplete gamma ratio function,
$\kappa_{i}={\cal{U}}_{j;\alpha}(x_{i}),i=1,2,j \in \{1,2,3\},{\cal{U}}_{1;\alpha}(x_{i})=x_{i}^{-\alpha},x_{i}\leq 0;{\cal{U}}_{2;\alpha}(x_{i})=(-x_{i})^{\alpha}, x_{i}>0,\alpha>0$ and ${\cal{U}}_3(x_{i})={\mathcal{U}}_{3;\alpha}(x_{i})=e^{-x_{i}},~
\forall\ x_{i}.$\\
{\bf Remark 2.1.} The continuity of
the limit df  $\hat{\Phi}_{r,s}^{(m,k)}(x_{}
,y_{})$ in $(i)$ implies the
continuity of the limit $\hat{\Psi}_{r,s}^{(m,k)}(x, y).$ Hence the
convergence in $(iii)$ is uniform with respect to  $x$ and $y.$\\
{\bf Remark 2.2.} It is natural to look for the limitations on
$\nu_n,$ under which we get the relation  $\hat{\Phi}_{r,s}^{(m,k)}(x_{},y_{})\equiv \hat{\Psi}_{r,s}^{(m,k)}(x, y)~\forall~x,y.$ In view of Theorem 2.1,
the last equation is satisfied if and only if  the df $H(z)$ is
degenerate at one, which means the asymptotically almost
randomlessness of $\nu_n.$ We assume, due to Remark 2.2, that $H(z)$ is a non-degenerate df  and $H(+0) = 0,$ i.e., continuous at zero.\\
{\bf Proof of the implication $(i)+(ii)\Rightarrow(iii)$:} First, we note that $~\Phi_{\grave{r},\grave{s}:n}^{(m,k)}(x_{n}
,y_{n})$ can be written in the form (see, Theorem 2.3 in Barakat et al.,  2014b)
$$\Phi_{\grave{r},\grave{s}:n}^{(m,k)}({x_{n}},{y_{n}})=1-\Gamma_{R_{r}}({N}{\ov{L}_{m}}(x_{n})) -\frac{1}{\Gamma(R_{r})}\int_{{N}{\ov{L}_{m}}(x_{n})}^{{N}}I_{\frac{{N}{\ov{L}_{m}}(y_{n})}{u}}(R_{s},R_{r}-R_{s})u^{R_{r}-1}{e^{-u}}du,\eqno (2.1)$$
 where ${\ov{L}_{m}}(.)=1-{{L}_{m}}(.).$ Now by using the total probability rule we get,
$$~\Phi_{\grave{r},\grave{s}:\nu_{n}}^{(m,k)}(x_{n}
,y_{n})=\sum_{t=r}^{\infty}\Phi_{\grave{r},\grave{s}:t}^{(m,k)}(x_{n}
,y_{n})P(\nu_{n}=t)
.\eqno (2.2)~$$  Assume that
$~H_n(z)=\sum_{t\leq z} P(\nu_{n}=t)=P(\nu_{n}\leq z)$  and  $z=[\frac{t}{n}],~$ where  $[\theta]$ denotes the greatest integer
part of $\theta.$ Thus, the relation
(2.1) show that the sum term in (2.2) is a Riemann sum of the integral
$$~\Phi_{\grave{r},\grave{s}:\nu_{n}}^{(m,k)}(x_{n}
,y_{n})=\int_{0}^{\infty}
\Phi_{\grave{r},\grave{s}:n}^{(m,k)}({x_{n}},{y_{n}},z)dH_{n}(nz),\eqno (2.3)~$$
where, for sufficiently large $n,$ we have $$\Phi_{\grave{r},\grave{s}:n}^{(m,k)}({x_{n}},{y_{n}},z)=1-\Gamma_{R_{r}}(z\grave{N}{\ov{L}_{m}}(x_{n})) -\frac{1}{\Gamma(R_{r})}\int_{z\grave{N}{\ov{L}_{m}}(x_{n})}^{z\grave{N}}I_{\frac{z\grave{N}{\ov{L}_{m}}(y_{n})}{u}}(R_{s},R_{r}-R_{s})u^{R_{r}-1}{e^{-u}}du,$$
where $\grave{N}=(\frac{\ell-1}{z}+n)\sim n.$ Appealing to the condition $(i),$ Theorem 2.3 in Barakat et al. (2014b) and Remark 2.2, we get
$$\Phi_{\grave{r},\grave{s}:n}^{(m,k)}({x_{n}},{y_{n}},z)
{\renewcommand{\arraystretch}{0.1}
\begin{array}{c}
\stackrel{\textstyle w}{\longrightarrow}\\
\scriptstyle{n}\end{array}}\overline{\Omega}_{r,s}^{(m,k)}(z\kappa_{1},z\kappa_{2}),\eqno (2.4)~$$ where the convergence is uniform with respect to $x$ and $y,$ over any finite interval of $z.$\\
Now, let $\xi$ be a continuity point of $H(z)$ such that
$1-H(\xi)<\epsilon,$ ($\epsilon$ is an aribtary small value). Then, we have
$$~\int_{\xi}^{\infty}\overline{\Omega}_{r,s}^{(m,k)}(z\kappa_{1},z\kappa_{2})dH(z)\leq \int_{\xi}^{\infty}dH(z)= 1-H(\xi)<\epsilon.\eqno (2.5)~$$
In view of the condition $(ii)$ we get, for sufficiently
large $n,$ that
$$~\int_{\xi}^{\infty}\Phi_{\grave{r},\grave{s}:n}^{(m,k)}({x_{n}},{y_{n}},z)
dH_n(nz)\leq 1-H_n (n\xi)\leq
2(1-H(\xi))<2\epsilon.\eqno (2.6)~$$  On the other hand, by the triangle inequality, we get
$$\left|\int_{0}^{\xi}\Phi_{\grave{r},\grave{s}:n}^{(m,k)}({x_{n}},{y_{n}},z)dH_n(nz)-\int_{0}^{\xi}
\overline{\Omega}_{r,s}^{(m,k)}(z\kappa_{1},z\kappa_{2})dH(z)\right|~$$ $$ \leq
\left|\int_{0}^{\xi}\Phi_{\grave{r},\grave{s}:n}^{(m,k)}({x_{n}},{y_{n}},z)dH_n(nz)-\int_{0}^{\xi}
\overline{\Omega}_{r,s}^{(m,k)}(z\kappa_{1},z\kappa_{2})dH_n(nz)\right| $$   $$+\left|\int_{0}^{\xi}
\overline{\Omega}_{r,s}^{(m,k)}(z\kappa_{1},z\kappa_{2})dH_n(nz)-\int_{0}^{\xi} \overline{\Omega}_{r,s}^{(m,k)}(z\kappa_{1},z\kappa_{2})dH(z)\right|,\eqno (2.7) $$ where the convergence in (2.4) is uniform
over the finite interval $0\leq z\leq \xi.$ Therefore, for arbitrary
$\epsilon>0$ and for sufficiently large $n,$ we have $$\left|\int_{0}^{\xi}\left[\Phi_{\grave{r},\grave{s}:n}^{(m,k)}({x_{n}},{y_{n}},z)-\overline{\Omega}_{r,s}^{(m,k)}(z\kappa_{1},z\kappa_{2})\right]dH_{n}(nz)\right|\leq\epsilon H_n(n\xi)\leq\epsilon.\eqno (2.8)~$$ In order to estimate the second difference on the right hand
side of (2.7), we construct Riemann sums which are close to the
integral there. Let $T$ be a fixed number and
$0=\xi_0<\xi_1<...<\xi_T=\xi$ be  continuity points of
$H(z).$ Furthermore, let $T$ and $\xi_{i}$ be such that
$$\left|\int_{0}^{\xi} \overline{\Omega}_{r,s}^{(m,k)}(z\kappa_{1},z\kappa_{2})dH_n(nz)-
\sum_{i=1}^{T}\overline{\Omega}_{r,s}^{(m,k)}(\xi_i \kappa_{1}
,\xi_i \kappa_{2}) (H_n(n\xi_i)-H_n(n\xi_{i-1}))\right|<\epsilon~$$
and
$$\left|\int_{0}^{\xi}\overline{\Omega}_{r,s}^{(m,k)}(z\kappa_{1},z\kappa_{2})dH(z)-
\sum_{i=1}^{T}\overline{\Omega}_{r,s}^{(m,k)}(\xi_i \kappa_{1},\xi_i \kappa_{2}) (H(\xi_i)-H(\xi_{i-1}))\right|<\epsilon.~$$
Since, by the assumption
$~H_n(n\xi_{i}){\renewcommand{\arraystretch}{0.1}
\begin{array}{c}
\stackrel{\textstyle w}{\longrightarrow}\\
\scriptstyle{n}\end{array}}H(\xi_{i}),~0\leq i\leq T,$ the two
Riemann sums are closer to each other than $\epsilon$ for all $n$
sufficiently large. Thus, once again by the triangle inequality,
the absolute value of the difference of the integrals is smaller
than $3\epsilon.$ Combining this fact with (2.8), the left hand
side of (2.7) becomes smaller than $4\epsilon$ for all large $n.$
Therefore, in view of (2.5), (2.6) and (2.4), we have
$$~\left| \Phi_{\grave{r},\grave{s}:\nu_n}^{(m,k)}(x_{n}
,y_{n})-\hat{\Psi}_{r,s}^{(m,k)}(x,y_{})\right|
 <\left|\int_{0}^{\xi}\Phi_{\grave{r},\grave{s}:n}^{(m,k)}({x_{n}},{y_{n}},z)dH_n(nz)-
 \int_{0}^{\xi} \overline{\Omega}_{r,s}^{(m,k)}(z\kappa_{1},z\kappa_{2})\\dH(z)\right|~$$
$$~+\int_{\xi}^{\infty}\Phi_{\grave{r},\grave{s}:n}^{(m,k)}({x_{n}},{y_{n}},z)dH_n(nz)+
\int_{\xi}^{\infty}\overline{\Omega}_{r,s}^{(m,k)}(z\kappa_{1},z\kappa_{2})dH(z)<7\epsilon.~$$
This completes the proof of the first part of the  theorem.\\
{\bf Proof of the implication $~(i)+(iii)\Rightarrow (ii)$:} Starting with the relation (2.3), we select a subsequence $\{n'\}$ of
$\{n\}$ for which $~H_{n'}(n'z)~$ converges weakly to an extended
df ${H'}(z)$ (i.e., ${H'}(\infty)-H'(0)\leq 1$ and such a subsequence
exists by the compactness of df's). Then, by repeating the first
part of the theorem for the subsequence $\{n'\},$ with the
exception that we choose $\xi$ so that
$H'(\infty)-H'(\xi)<\epsilon,$ we get $\hat{\Psi}_{r,s}^{(m,k)}(x,y)=\int_{0}^{\infty}\overline{\Omega}_{r,s}^{(m,k)}(z\kappa_{1},z\kappa_{2})dH'(z).$ Since
the function $\hat{\Psi}_{r,s}^{(m,k)}(x,y)$ is a df, we get
$\hat{\Psi}_{r,s}^{(m,k)}(\infty,\infty)=1=\int_{0}^{\infty}
dH'(z)=H'(\infty)-H'(0),$ which implies that $H'(z)$ is a df.
Now, if $H_n(n z)$ did not converge weakly, then we can select two
subsequences $\{n'\}$ and $\{n''\}$ such that
$H_{n'}(n'z){\renewcommand{\arraystretch}{0.1}
\begin{array}{c}
\stackrel{\textstyle w}{\longrightarrow}\\
\scriptstyle{n'}\end{array}}H'(z)$ and
$H_{n''}(n''z){\renewcommand{\arraystretch}{0.1}
\begin{array}{c}
\stackrel{\textstyle w}{\longrightarrow}\\
\scriptstyle{n''}\end{array}}H''(z),$ where $H'(z)$ and $H''(z)$
are df's. In this case, we get
$$\hat{\Psi}_{r,s}^{(m,k)}(x,y_{})=\int_{0}^{\infty}\overline{\Omega}_{r,s}^{(m,k)}(z\kappa_{1},z\kappa_{2})dH'(z)\!\!=\!\!
\int_{0}^{\infty}\overline{\Omega}_{r,s}^{(m,k)}(z\kappa_{1},z\kappa_{2})dH''(z).$$ Thus, let ($ y\rightarrow \infty $), we get
$$~\int_{0}^{\infty}\Gamma_{R_{r}}(z \kappa_{1}^{m+1})dH'(z) =
\int_{0}^{\infty}\Gamma_{R_{r}}(z \kappa_{1}^{m+1})dH''(z).\eqno (2.9)~$$
Appealing to equation $ (2.9)$ and by using the same argument which is applied in the proof of the second part of Theorem 2.1 in Barakat, 1997,
 we can easily prove $H'(z)=H''(z).$ This complete the proof of the second part.\\
{\bf Proof of the implication $(ii)+(iii)\Rightarrow(i)$:}  For proving this part, we need first present the following lemma.\\
{\bf Lemma 2.1.} For all $x_{i},i=1,2,$ we have
$$[1-\Gamma_{R_{r_{i}}}(N \ov{L}_{m}(x_{in}))]-\sigma_{i,N} \leq P({Z^{(n)}_{\grave{r}_{i}:n}}< x_{i}) \leq [1-\Gamma_{R_{r_{i}}}(N\ov{L}_{m}(x_{in}))]+\rho_{i,N}\eqno(2.10)$$
and $$\Gamma_{R_{r_{i}}}(N\ov{L}_{m}(x_{in}))-\rho_{i,N}\leq P({Z^{(n)}_{\grave{r}_{i}:n}}\geq x_{i}) \leq \Gamma_{R_{r_{i}}}(N\ov{L}_{m}(x_{in}))+\sigma_{i,N},\eqno(2.11)$$
where $0<\rho_{i,N} ~,\sigma_{i,N}~ {\renewcommand{\arraystretch}{0.1}
\begin{array}{c}
\stackrel{\textstyle }{\longrightarrow}\\
\scriptstyle{N}\end{array}} ~ 0$  (or equivalently, as $n\longrightarrow\infty),(R_{r_{1}},R_{r_{2}})=(R_{r_{}},R_{s})$ and $(r_{1},r_{2})=(r,s).$\\
{\bf Proof.} Since $0 \leq \Gamma_{R_{r_{i}}}(x)\leq 1 ~\forall~ x,$ the proof of the lemma  will immediately follow from
the result of Smirnov (1952) (Theorem 3, p. 133, or Lemma   2.1 in Barakat, 1997).

We now turn to the proof of the last part of Theorem 2.1. In view of Remark 2.1, we can assume,
without any loss of generality, that the df
$ \hat{\Psi}_{{r},{s}}^{(m,k)}(x_{}
,y_{})$ is continuous. Therefore, the condition $(iii)$ will be
satisfied for all  univariate marginals of $ \hat{\Psi}_{{r},{s}}^{(m,k)}(x_{}
,y_{}),$
i.e., we have
$$~\Phi_{\grave{r_{i}}:\nu_n} ^{(m,k)}{(x_{in})}{\renewcommand{\arraystretch}{0.1}
\begin{array}{c}
\stackrel{\textstyle w}{\longrightarrow}\\
\scriptstyle{n}\end{array}}\hat{\Psi}_{{r_{i}}} ^{(m,k)}{(x_{i})},i=1,2,\eqno
(2.12)~$$ where $\hat{\Psi}_{{r_{i}}} ^{(m,k)}{(x_{i})}=\int_{0}^{\infty}[1-\Gamma_{R_{r_{i}}}(z\kappa_{i}^{m+1})]dH(z),i=1,2,$ is the  marginals
df's of $\hat{\Psi}_{{r},{s}}^{(m,k)}(x_{}
,y_{}).$ We shall now prove
$$\Phi_{\grave{r_{i}}:n} ^{(m,k)}{(x_{in})}{\renewcommand{\arraystretch}{0.1}
\begin{array}{c}
\stackrel{\textstyle w}{\longrightarrow}\\
\scriptstyle{n}\end{array}} \hat{\Phi}_{{r_{i}}} ^{(m,k)}{(x_{i})}=1-\Gamma_{R_{r_{i}}}(\kappa_{i}^{m+1}),i=1,2. \eqno (2.13)~$$
In view of Lemma 2.1, we first show that the sequence
$\{Z^{(n)}_{\grave{r_{i}}:n}\}_n,~i=1,2,$ is
stochastically bounded (see, Feller, 1979). If we assume the contrary, we would find
$\varepsilon_{i,1},\varepsilon_{i,2}>0$ such that at least one of
the two following relations \\ (a) $\qquad {\renewcommand{\arraystretch}{0.1}
\begin{array}{c}
~\\
~\\
~\\
{\overline {\lim}}\\
\scriptstyle{n\to\infty}\end{array}}~P(Z^{\langle
n\rangle}_{\grave{r_{i}}:n}\geq x_{i}) \geq\varepsilon_{i,1}>0, ~~\forall~~
x_{i}>0,~ i=1,2,~$\\ (b) $\qquad
{\renewcommand{\arraystretch}{0.1}
\begin{array}{c}
~\\
~\\
~\\
{\overline {\lim}}\\
\scriptstyle{n\to\infty}\end{array}}P(Z^{\langle
n\rangle}_{\grave{r_{i}}:n}<x_{i}) \geq\varepsilon_{i,2}>0, ~~\forall
~~x_{i}<0,~ i=1,2~$\\ is satisfied. The assertions (a) and (b)
mean that the sequence $\{Z^{(n)}_{\grave{r_{i}}:n}\}_n,$  is not stochastically bounded
at the left $(-\infty)$ and at the right $(+\infty),$
respectively. Let the assumption (a) be true. Since $H(z)$ is
non-degenerate df, we find $\varepsilon_0>0$ and $\beta>0$ such
that
$$~P\left(\frac{\nu_n}{n}\geq\beta\right)\geq\varepsilon_0,~
\mbox{for sufficiently large}~n.\eqno (2.14)~$$ Using the following well
known inequality, for $i=1,2,$
$$~P\left(Z_{\grave{r_{i}}:n_{1}}^{(n)}\geq x_i\right)\geq P\left(Z_{\grave{r_{i}}:n_{2}}^{\langle
n\rangle}\geq x_i \right ),\quad \forall ~n_{1} \geq n_{2}.\eqno
(2.15)~$$ We thus get  the following inequalities, for
sufficiently large $n,$
$$P\left(Z_{\grave{r_{i}}:\nu_n}^{(n)}\geq x_i\right)
\geq \sum_{t\geq[n\beta]}P\left(Z_{\grave{r_{i}}:t}^{(n)}\geq
x_i \right )P(\nu_n=t)~$$ $$~\geq P\left(Z_{\grave{r_{i}}:[n\beta]}^{\langle
n\rangle}\geq x_i \right )P\left(\nu_n\geq [n\beta] \right )\geq
\varepsilon_0 P\left(Z_{\grave{r_{i}}:[n\beta]}^{(n)}\geq x_i
\right ),~i=1,2,$$ (note that $P\left(\nu_n\geq
[n\beta]\right)\geq P\left(\nu_n\geq n\beta \right ) $).
Therefore,
$${\renewcommand{\arraystretch}{0.1}
\begin{array}{c}
~\\
~\\
~\\
{\overline {\lim}}\\
\scriptstyle{n\to\infty}\end{array}}P(Z^{\langle
n\rangle}_{\grave{r_{i}}:\nu_n}\geq x_i)
\geq\varepsilon_0{\renewcommand{\arraystretch}{0.1}
\begin{array}{c}
~\\
~\\
~\\
{\overline {\lim}}\\
\scriptstyle{n\to\infty}\end{array}} P(Z^{(n)}_{\grave{r_{i}}:[n\beta]}\geq
x_i).$$  Now, if we find $\varepsilon'_{i,1}>0$ such that $
{\renewcommand{\arraystretch}{0.1}
\begin{array}{c}
~\\
~\\
~\\
{\overline {\lim}}\\
\scriptstyle{n\to\infty}\end{array}} P(Z^{\langle
n\rangle}_{\grave{r_{i}}:[n\beta]}\geq x_i)\geq \varepsilon'_{i,1}>0,$ we get
${\renewcommand{\arraystretch}{0.1}
\begin{array}{c}
~\\
~\\
~\\
{\overline {\lim}}\\
\scriptstyle{n\to\infty}\end{array}} P(Z^{\langle
n\rangle}_{\grave{r_{i}}:\nu_n}\geq x_i)\geq\epsilon_0 \varepsilon'_{i,1}>0,$
which contradicts the right stochastic boundedness of the sequence
$\{Z_{\grave{r_{i}}:\nu_n}^{(n)}\}_n$ and consequently
contradicts the relation (2.12). However, if such an
$\varepsilon'_{i,1}>0$ does not exist we have
${\renewcommand{\arraystretch}{0.1}
\begin{array}{c}
~\\
~\\
~\\
{\overline {\lim}}\\
\scriptstyle{n\to\infty}\end{array}} P(Z^{\langle
n\rangle}_{\grave{r_{i}}:[n\beta]}\geq x_i)=0,$ which in view of Lemma 2.1 (relation (2.10)) leads to the
following chain of implications ($\forall~x_i>0$) $P(Z^{\langle
n\rangle}_{\grave{r_{i}}:[n\beta]}\geq x_i)\rightarrow 0 \Rightarrow
\Gamma_{R_{r_{i}}}([N\beta]\ov{L}_{m}(x_{in}))\rightarrow 0 \Rightarrow
[N\beta]\ov{L}_{m}(x_{in}) \rightarrow 0\Rightarrow
N\ov{L}_{m}(x_{in})\rightarrow 0$ (since $N\ov{L}_{m}(x_{in})\rightarrow 0)$
$\Rightarrow \Gamma_{R_{r_{i}}}(N\ov{L}_{m}(x_{in}))\rightarrow 0 \Rightarrow
P(Z^{(n)}_{\grave{r_{i}}:n}\geq x_i)\rightarrow 0,$ which
contradicts the assumption (a). Consider the assumption (b). Since
$H(z)$ is a df we can find a positive integer $\delta$ and real
number $\alpha >0$ such that
$$~P(\frac {\nu_n}{n} \leq \delta)\geq\alpha,\qquad\mbox{ for sufficiently
large}~ n.\eqno (2.16)~$$ Therefore, in view of (2.16) and the
inequality (2.15), we have
$$P(Z_{\grave{r_{i}}:\nu_n}^{(n)}<x_i)\geq\sum_{t=r}^{\delta n}
P(Z_{\grave{r_{i}}:t}^{(n)}<x_i)P(\nu_n=t)\geq P(Z_{\grave{r_{i}}:\delta
n}^{(n)}<x_i)P(\frac {\nu_n}{n} \leq \delta)~$$
$$~\geq \alpha P(Z_{\grave{r_{i}}:\delta n}^{(n)}<x_i),~
i=1,2.~$$ Hence, we get ${\renewcommand{\arraystretch}{0.1}
\begin{array}{c}
~\\
~\\
~\\
{\overline {\lim}}\\
\scriptstyle{n\to\infty}\end{array}} P(Z^{\langle
n\rangle}_{\grave{r_{i}}:\nu_n}< x_i)
\geq\alpha{\renewcommand{\arraystretch}{0.1}
\begin{array}{c}
~\\
~\\
~\\
{\overline {\lim}}\\
\scriptstyle{n\to\infty}\end{array}} P(Z^{\langle
n\rangle}_{\grave{r_{i}}:\delta n}< x_i).~$ By using  Lemma 2.1 (relation (2.11)) and applying the same
argument as in the case (a), it is easy to show that the last
inequality leads to a contradiction (the last inequality, in view
of the assumption (b)), which yields that the sequences
$\{Z_{\grave{r_{i}}:{n}}^{(n)}\}_n,i=1,2,~$ is not
stochastically bounded at the left. This completes the proof that
the sequences $~\{Z_{\grave{r_{i}}:n}^{(n)}\}_n,~i=1,2,~$
are stochastically bounded. Now, if  $\Phi_{\grave{r_{i}}:n} ^{(m,k)}{(x_{in})},i=1,2,$ did not converge weakly, then we could select two subsequences $\{n'\}$ and $\{n''\}$ such that $\Phi_{\grave{r_{i}}:n'} ^{(m,k)}{(x_{in'})}$ would converge weakly to $~\hat{\Phi}_{{r_{i}}} ^{'(m,k)}{(x_{i})}~$ and
$~\Phi_{\grave{r_{i}}:n''} ^{(m,k)}{(x_{in''})}~$ to another limit df
$\hat{\Phi}_{{r_{i}}} ^{''(m,k)}{(x_{i})}.$ In this case we get (by
repeating the first part of Theorem 2.1 for the univariate case
and for the two subsequences $\{n^{'}\},\{n^{''}\}$)
$$~{\overline{\Omega}}^{({m},k)}_{r_{i}}(x_{i})=\int_{0}^{\infty}\left[1-\Gamma_{R_{r_{i}}}(z\kappa_{i}^{'(m+1)})\right]dH(z)=\int_{0}^{\infty}\left[1-\Gamma_{R_{r_{i}}}(z\kappa_{i}^{''(m+1)})\right]dH(z).~$$
However, Lemma 3.2 in Barakat (1997) shows that the last
equalities, cannot hold unless $\kappa_{i}^{'}\equiv \kappa_{i}^{''},i=1,2.$ Hence
the relation (2.13) is proved. Hence, the proof of Theorem
2.1 is completed.~\rule {2mm}{2mm}

Let ${\cal G}^{}$ and ${\cal G}_{\nu}^{}$ be the classes of all possible limit df's in $(i)$ and $(iii),$ respectively. The class ${\cal G}^{}$ is fully determined
by Barakat et al. (2014b). Furthermore, let $S^{}$ and $S_{\nu}^{}$ be the
necessary and sufficient conditions for the validity of the relations  $(i)$ and $(iii),$ respectively. The
following corollary characterizes the class ${\cal G}_\nu^{}.$\\
{\bf Corollary 2.1.} For every df  $\hat{\Psi}^{(m,k)}(x,y)$ in
${\cal G}_\nu^{}$  there exists a unique df $\hat{\Phi}^{(m,k)}(x,y)$ in
${\cal G}^{},$  such that $\hat{\Psi}^{(m,k)}(x,y)$
is uniquely determined by the representation $(iii).$ Moreover, $S^{}=S_{\nu}^{}.$\\
{\bf Proof of corollary 2.1.} Let us first prove the implication $\{\hat{\Phi}_{r,s}^{'(m,k)}(x,y)\neq \hat{\Phi}_{r,s}^{''(m,k)}(x,y)\}\\\Longrightarrow \{{\hat{\Psi}}_{r,s}^{'(m,k)}(x,y)\neq {\hat{\Psi}}_{r,s}^{''(m,k)}(x,y)\}.$
If we assume the contrary, we get $ {\hat{\Psi}}_{r,s}^{'(m,k)}(x,y)= {\hat{\Psi}}_{r,s}^{''(m,k)}(x,y),$
while $\hat{\Phi}_{r,s}^{'(m,k)}(x,y)\neq \hat{\Phi}_{r,s}^{''(m,k)}(x,y).$
 Appealing to the first part of Theorem 2.1 , we get
$\int_{0}^{\infty}\left[1-\Gamma_{R_{r_{i}}}(z\kappa_{i}^{'(m+1)})\right]dH(z)=\int_{0}^{\infty}\left[1-\Gamma_{R_{r_{i}}}(z\kappa_{i}^{''(m+1)})\right]dH(z),i=1,2.$ The last equalities,
as we have seen before, from Lemma 3.2 in Barakat (1997), cannot hold unless
$\kappa'_i=\kappa''_i,i=1,2.$ Therefore, Corollary 2.1 is followed as a consequence of
Theorem 2.1 and the last implication. This completes the proof of Corollary 2.1.~\rule {2mm}{2mm}\\
{\bf Theorem 2.2.} Consider the following three conditions :
 $$ \Phi_{r,s:n}^{(m,k)}(x_{n}
,y_{n})=P( Z^{(n)}_{{r},{s}:{n} }<\mathbf{x})=P(Z_{{r}:n}^{ (n)}<x_{},Z_{{{s}}:n}^{ (n)}<y){\renewcommand{\arraystretch}{0,1} \begin{array}{c}
\stackrel{w}{\longrightarrow}\\
\scriptstyle{n}\end{array}}\Phi_{r,s}^{(m,k)}(x_{}
,y_{}),x\leq y, \eqno (i)~$$
$$~ H_n(nz)= P(\frac{\nu_n}{n}<z){\renewcommand{\arraystretch}{0.1}
\begin{array}{c}
\stackrel{\textstyle w}{\longrightarrow}\\
\scriptstyle{n}\end{array}}H(z),~\eqno (ii)~$$
$$ \Phi_{r,s:\nu_n}^{(m,k)}(x_{n}
,y_{n})=P( Z^{(n)}_{{r},{s}:\nu_{n} }<\mathbf{x})=P(Z_{{r}:\nu_{n}}^{ (n)}<x_{},Z_{{{s}}:\nu_{n}}^{ (n)}<y){\renewcommand{\arraystretch}{0,1} \begin{array}{c}
\stackrel{w}{\longrightarrow}\\
\scriptstyle{n}\end{array}}
    \Psi_{r,s}^{(m,k)}(x_{}
,y_{})$$
$$=\int_{0}^{\infty} \underline{\Omega}_{r,s}^{(m,k)}(z\rho_{1}
,z\rho_{2})dH(z).~\eqno (iii)~$$
 Then any two of the above conditions imply the remaining one, where $x_{n}=c_n {x_{}}+ d_n, y_{n}=c_n {y_{}}+ d_n,c_n>0, d_n$ are suitable normalizing constants,  $1\leq r<s\leq n,Z_{{r_{i}}:n}^{(n)}=\frac{X^{}_{{r_{i}}:n}-d_n}{c_n}, i=1,2, \Phi_{r,s}^{(m,k)}(x_{}
,y_{})$ is a non-degenerate df, $~H(z)~$
is a df with $~H(+0)=0,$
$$\underline{\Omega}_{r,s}^{(m,k)}(\rho_{1}
,\rho_{2})=\left\{
    \begin{array}{ll}
    \Gamma_{{s}}( \rho_{2}),& x\geq y,\\ \frac{1}{(r-1)!} \int_{0}^{\rho_{1}}\Gamma_{{s}-{r}}(\rho_{2}-u)u^{r-1}{e^{-u}}du,& x\leq y,
     \end{array}
      \right.$$
$\rho_{i}={\cal{V}}_{j;\beta}(x_{i}),i=1,2,j \in \{1,2,3\},{\cal{V}}_{1;\beta}(x_{i})=(-x_{i})^{-\beta},x_{i}\leq 0;{\cal{V}}_{2;\beta}(x_{i})=x_{i}^{\beta}, x_{i}>0,\beta>0$ and ${\cal{V}}_3(x_{i})={\mathcal{V}}_{3;\beta}(x_{i})=e^{x_{i}},~
\forall\ x_{i}.$\\
{\bf Theorem 2.3.} Consider the following three conditions :
   $$ \Phi_{{r},\grave{s}:n}^{(m,k)}(x_{n}
,y_{n})=P( Z^{(n)}_{{r},\grave{s}:{n} }<\mathbf{x})=P(Z_{{r}:n}^{ (n)}<x_{},Z_{{{\grave{s}}}:n}^{ (n)}<y){\renewcommand{\arraystretch}{0,1} \begin{array}{c}
\stackrel{w}{\longrightarrow}\\
\scriptstyle{n}\end{array}}{\Phi}_{r}^{(m,k)}{(x)}{\hat{\Phi}_{{s}}^{(m,k)}(y)}$$
$$ ~~~~= \Gamma_{{r}}(\rho_{1}))[1-\Gamma_{R_{s}}(\kappa_{2}^{m+1}),1\leq r,s\leq n, \eqno (i)$$
$$~ H_n(nz)= P(\frac{\nu_n}{n}<z){\renewcommand{\arraystretch}{0.1}
\begin{array}{c}
\stackrel{\textstyle w}{\longrightarrow}\\
\scriptstyle{n}\end{array}}H(z),~\eqno (ii)~$$
$$ \Phi_{{r},\grave{s}:\nu_n}^{(m,k)}(x_{n}
,y_{n})=P( Z^{(n)}_{{r},\grave{s}:\nu_{n} }<\mathbf{x})=P(Z_{{r}:\nu_{n}}^{ (n)}<x_{},Z_{{{\grave{s}}}:\nu_{n}}^{ (n)}<y){\renewcommand{\arraystretch}{0,1} \begin{array}{c}
\stackrel{w}{\longrightarrow}\\
\scriptstyle{n}\end{array}}
    \underline{\Omega}_{r}^{(m,k)}{(z\rho_{1})}{\overline{\Omega}_{{s}}^{(m,k)}(z\kappa_{2})}.~\eqno (iii)~$$
Then any two of the above conditions imply the remaining one,  where $~x_{n}=c_n {x}+ d_n,y_{n}=a_n {y}+ b_n, a_n,c_n>0, b_n, d_n$ are suitable normalizing constants, ${ Z}^{(n)}_{ {r}:n}=\frac{X^{}_{{r}:n}-d_n}{c_n},{ Z}^{(n)}_{{\grave{s}}:n}=\frac{X^{}_{{\grave{s}}:n}-b_n}{a_n},
 {\Phi}_{r}^{(m,k)}{(x)},{\hat{\Phi}_{{s}}^{(m,k)}(y)}~$ are non-degenerate df's,  $~H(z)~$
is a df with $~H(+0)=0,$
$$\underline{\Omega}_{r}^{(m,k)}{(z\rho_{1})}=\int_{0}^{\infty} \Gamma_{{r}}(z\rho_{1})dH(z) ~\mbox{and}~
{\overline{\Omega}_{{s}}^{(m,k)}(z\kappa_{2})}=\int_{0}^{\infty}\left[1-\Gamma_{R_{s}}(z\kappa_{2}^{m+1})\right]dH(z).$$
{\bf Proof of Theorems 2.2 and  2.3.} Without significant modifications, the method of the proof
of Theorems 2.2 and 2.3 are the same as that  Theorem
2.1, except only the obvious changes. Hence, for brevity the details of the proof are omitted.~\rule {2mm}{2mm}
\subsection{ The interrelation of $\nu_n$ and the basic rv's is not
restricted}
When the interrelation between the random index and the basic
variables is not restricted, parallel theorem of Theorem 2.1
may be proved by  replacing the condition $(ii)$ by a stronger one. Namely, the weak convergence of
the df $H_n(nz)$  must
be replaced by the convergence in probability of
the rv $\frac{\nu_n}{n}$  to a positive rv $\mathcal{T}.$  However, the key ingredient of the proof of this parallel result
is to prove the mixing property, due to
R\'{e}nyi (see, Barakat and Nigm, 1996) of the sequence
of order statistics under consideration. In the sense of R\'{e}nyi
a sequence $\{{u}_n\}$ of rv's is  called mixing if for any event ${ \mathcal{E}}$ of positive
probability, the conditional df of $\{u_n\},$ under the
condition ${\mathcal{ E}},$ converges  weakly to a non-degenerate df, which
does not depend on ${ \mathcal{E}},$ as $n\to\infty.$ The following lemma
proves the mixing property for the sequence $\{Z_{\grave{r},\grave{s}:n}^{(n)}\}_n. $
\\\textbf{Lemma 2.2.} Under the condition $(i)$ in Theorem 2.1 the sequence $\{Z_{\grave{r},\grave{s}:n}^{(n)}\}_n $ is mixing.
\\\textbf{Proof.} The lemma will be proved if one shows the relation $ P(Z_{\grave{r},\grave{s}:n}^{(n)}< \mathbf{x} \mid Z_{\grave{r},\grave{s}:l}^{(l)}<\mathbf{x} ){\renewcommand{\arraystretch}{0.1}
\begin{array}{c}
\stackrel{\textstyle w}{\longrightarrow}\\
\scriptstyle{n}\end{array}}\hat{\Phi}_{{r},{s}}^{(m,k)}({x},y),$ for all integers $l= r,r+1,....$  The sufficiency of the above relation can easily be proved as a direct multivariate extension of Lemma 6.2.1, of Galambos (1987). However, this relation is equivalent to $$P(Z_{\grave{r},\grave{s}:n}^{(n)}\geq \mathbf{x}\mid Z_{\grave{r},\grave{s}:l}^{\langle l \rangle}\geq \mathbf{x} ){\renewcommand{\arraystretch}{0.1}
\begin{array}{c}
\stackrel{\textstyle w}{\longrightarrow}\\
\scriptstyle{n}\end{array}}\ov{\Phi}_{{r},{s}}^{(m,k)}({x},y),\eqno (2.17)~$$ where
$\ov{\Phi}_{{r},{s}}^{(m,k)}({x},y)$ is the survival function of the limit df $\hat{\Phi}_{{r},{s}}^{(m,k)}({x},y),$ i.e.,
$$~\ov{\Phi}_{{r},{s}}^{(m,k)}({x},y)= 1-\hat{\Phi}_{{r}}^{(m,k)}({x})-\hat{\Phi}_{{s}}^{(m,k)}(y)+\hat{\Phi}_{{r},{s}}^{(m,k)}({x},y).$$
Therefore, our lemma will be established if one proves the relation (2.17). Now, we can write $$ P(Z_{\grave{r},\grave{s}:n}^{(n)}\geq \mathbf{x}\mid Z_{\grave{r},\grave{s}:l}^{(l)}\geq \mathbf{x})=P(Z_{\grave{r},\grave{s}:n}^{(n)}\geq \mathbf{x},Z_{\grave{r},\grave{s}:l}^{(n)}<\mathbf{x}\mid Z_{\grave{r},\grave{s}:l}^{(l)}\geq \mathbf{x})~$$ $$~+P(Z_{\grave{r},\grave{s}:n}^{(n)}\geq \mathbf{x},Z_{\grave{r},\grave{s}:l}^{(n)}\geq \mathbf{x}\mid Z_{\grave{r},\grave{s}:l}^{(l)}\geq \mathbf{x}).\eqno (2.18)~$$ Bearing in mind that all $X_{r:n},r=1,2,...,n,$ are i.i.d rv's, the first term in (2.18) can be written in the form $$ P(Z_{\grave{r},\grave{s}:n}^{(n)}\geq \mathbf{x},Z_{\grave{r},\grave{s}:l}^{(n)}<\mathbf{x}\mid Z_{\grave{r},\grave{s}:l}^{(l)}\geq \mathbf{x})=P(Z_{\grave{r},\grave{s}:(n-l)}^{*(n)}\geq \mathbf{x},Z_{\grave{r},\grave{s}:l}^{(n)}< \mathbf{x}\mid Z_{\grave{r},\grave{s}:l}^{(l)}\geq \mathbf{x})$$ $$=P(Z_{\grave{r},\grave{s}:(n-l)}^{*(n)}\geq \mathbf{x} )-P(Z_{\grave{r},\grave{s}:(n-l)}^{*(n)}\geq \mathbf{x},Z_{\grave{r},\grave{s}:l}^{(n)}\geq \mathbf{x}\mid Z_{\grave{r},\grave{s}:l}^{(l)}\geq \mathbf{x}),$$ where $$ Z_{\grave{r},\grave{s}:(n-l)}^{*(n)}=(Z_{\grave{r}:(n-l)}^{*(n)},Z_{\grave{s}:(n-l)}^{*(n)}),$$ $$Z_{\grave{r}:(n-l)}^{*(n)}=(({r}\mbox{th largest of}~X_{1,l+1:n-l},X_{1,l+2:n-l},...,X_{1,n:n-l})-b_{n})/a_{n},$$  $$Z_{\grave{s}:(n-l)}^{*(n)}=(({s}\mbox{th largest of}~ X_{2,l+1:n-l},X_{2,l+2:n-l},...,X_{2,n:n-l})-b_{n})/a_{n}.$$ Therefore, in view of (2.18), we have $$P(Z_{\grave{r},\grave{s}:n}^{(n)}\geq \mathbf{x}\mid Z_{\grave{r},\grave{s}:l}^{(l)}\geq \mathbf{x})= P(Z_{\grave{r},\grave{s}:(n-l)}^{*(n)}\geq \mathbf{x} )-\Delta_n(\mathbf{x}),\eqno (2.19)~$$ where $\Delta_n(\mathbf{x})= P(Z_{\grave{r},\grave{s}:n}^{(n)}\geq \mathbf{x},Z_{\grave{r},\grave{s}:l}^{(l)}\geq\mathbf{x}\mid Z_{\grave{r},\grave{s}:l}^{(l)}\geq \mathbf{x} )- P(Z_{\grave{r},\grave{s}:(n-l)}^{*(n)}\geq \mathbf{x}, Z_{\grave{r},\grave{s}:l}^{(n)}\geq \mathbf{x}\mid Z_{\grave{r},\grave{s}:l}^{(l)}\geq \mathbf{x} ).$ By using the well-known inequalities $~Z_{\grave{r},\grave{s}:(n-l)}^{*(n)} \leq Z_{\grave{r},\grave{s}:n}^{(n)}~$ and $~P(B \bigcap C)-P( A \bigcap C ) \leq P(B)-P(A),$ for any three events $ A, B~ \mbox{and}~C,$ for which $A\subseteq B,$ we get $$ 0\leq \Delta_n(\mathbf{x})P(Z_{\grave{r},\grave{s}:l}^{(l)}\geq\mathbf{x}) \leq P(Z_{\grave{r},\grave{s}:n}^{(n)}\geq \mathbf{x})- P(Z_{\grave{r},\grave{s}:(n-l)}^{*(n)}\geq \mathbf{x}).\eqno (2.20)~$$ On the other hand, by virtue of the condition $(i)$ in Theorem 2.1, it is easy to prove that $$ \lim_{n\rightarrow \infty }P(Z_{\grave{r},\grave{s}:(n-l)}^{*(n)}\geq \mathbf{x})= \lim_{n\rightarrow \infty }P(Z_{\grave{r},\grave{s}:(n-l)}^{(n)}\geq \mathbf{x})= \ov{\Phi}_{{r},{s}}^{(m,k)}({x},y)\eqno (2.21)~$$
(note that $N \ov{L}_{m}(x_{in}){\renewcommand{\arraystretch}{0.1}
\begin{array}{c}
\stackrel{\textstyle }{\longrightarrow}\\
\scriptstyle{n}\end{array}} \kappa_{i}^{m+1}\Rightarrow (N-l)\ov{L}_{m}(x_{in}){\renewcommand{\arraystretch}{0.1}
\begin{array}{c}
\stackrel{\textstyle }{\longrightarrow}\\
\scriptstyle{n}\end{array}} \kappa_{i}^{m+1},~\forall ~x_i's~$ for which $\kappa_i<\infty,~i= 1, 2~$). By combining the relations (2.19)-(2.21), the proof of the relation (2.17) follows immediately. Hence the required result.~\rule {2mm}{2mm}

Considering the facts that the normalizing constants, which may be used in the bivariate extreme case are the same as those for the univariate case, and the limit df $\hat{\Phi}_{{r},{s}}^{(m,k)}(x_{}
,y_{}) $  is continuous, we can easily by using Lemma 2.2,  show that the proof of the following theorem follows without any essential modifications as a direct multivariate extension of the proof of Theorem 2.1 in Barakat and El Shandidy (1990), except only the obvious changes.\\
\textbf{Theorem 2.4.} Consider the condition $$\qquad\qquad~~\frac{\nu_n}{n} {\renewcommand{\arraystretch}{0.1}
\begin{array}{c}
\stackrel{\textstyle p}{\longrightarrow}\\
\scriptstyle{n}\end{array}} ~\mathcal{T},~ \eqno (ii)'~$$
 where $\mathcal{T}~$ is a positive rv. Under the conditions of Theorem 2.1 , we have the implication $$(i)+(ii)'\Rightarrow (iii).$$
\section{Illustrative examples}
The range and midrange are widely used, particularly in statistical quality control as an
estimator of the dispersion tendency and in setting confidence intervals
for the population standard deviation as well as in Monte Carlo methods. In fact, the range itself is a very simple measure of dispersion,
gives a quick and easy to estimate indication about the spread of data. Let us defines the random generalized ranges  ${ \mathcal{R}}_{\nu_{n}}(m,k)\!=\!X^{}_{{\nu_n}:\nu_n}\!-\!X^{}_{1:\nu_n}$  and the random generalized midranges
${\mathcal{V}}_{\nu_{n}}(m,k)=\frac{X^{}_{1:\nu_n}+X^{}_{{\nu_n}:\nu_n}}{2},X_{1:\nu_n}= X(1,\nu_n,m,k)$ and $X_{\nu_{n}:\nu_{n}}=X(\nu_{n},\nu_{n},m,k).$
The normalized  generalized ranges and  the normalized
generalized midranges are defined by ${\cal R}^{(n)}_{\nu_{n}}(m,k)=A^{-1}_{n:r}({\cal R}_{\nu_{n}}(m,k)-B_{n:r})$ and ${\cal V}^{(n)}_{\nu_{n}}(m,k)=A^{-1}_{n:\upsilon}({\cal V}^{}_{\nu_{n}}(m,k)-B_{n:\upsilon}),$
respectively, where $A_{n:r}=2A_{n:\upsilon}=a_{n}> 0,B_{n:r}=b_{n}-d_{n}$ and $B_{n:\upsilon}=\frac{1}{2}(b_{n}+d_{n})$
are suitable normalizing constants. In this section, some illustrative examples for the most important distribution functions are
obtained, which lend further support to our theoretical results. In the following examples
we consider an important practical situation when $\nu_n$ has a geometric distribution with
mean $n.$ In this case we can easily show that $P(\frac{\nu_{n}}{n}< z){\renewcommand{\arraystretch}{0.1}
\begin{array}{c}
\stackrel{\textstyle w}{\longrightarrow}\\
\scriptstyle{n}\end{array}}H(z)=1-e^{-z}(z\geq 0).$\\
{\bf Example 3.1 (standard Cauchy distribution).}  Let $F(x)=\frac{1}{2}+ \frac{1}{\pi}\tan^{-1}{x}.$ Then $\alpha=\beta=1.$ In view of Theorems 1.1 and 2.1, Part 1, in Barakat et al. (2015b), we get, after some algebra,
$$\eta=\lim _{n \to \infty}\eta_{n}=\lim _{n \to \infty}\frac{a_{n}}{c_{n}} =~\left\{
    \begin{array}{ll} 1, ~~~~~~\mbox{if}~~ m=0,
\\0,~~~~~~ \mbox{if} ~~m>0,
\\\infty, ~~~~~\mbox{if}~ -1<m<0.\end{array}\right.$$
The random generalized ranges  and  midranges, for standard Cauchy distribution are given by, if $m=0,$
$$P({\cal R}^{(n)}_{\nu_{n}}(0,k)\leq r){\renewcommand{\arraystretch}{0,1} \begin{array}{c}
\stackrel{w}{\longrightarrow}\\
\scriptstyle{n}\end{array}}\int_{0}^{\infty}\int_{0}^{\infty}\left[1-\Gamma_{k}({z}({r-y})^{-1})\right] {y^{-2}}ze^{-z(1+y^{-1})}dydz$$
and  $$P({\cal V}^{(n)}_{\nu_{n}}(0,k)\leq v){\renewcommand{\arraystretch}{0,1} \begin{array}{c}
\stackrel{w}{\longrightarrow}\\\scriptstyle{n}\end{array}}\int_{0}^{\infty}\int_{-\infty}^{\infty}\left[1-\Gamma_{k}({z}({v-y})^{-1})\right] {y^{-2}}ze^{-z(1-y^{-1})}dydz,$$
 respectively, with $A_{n:r} = 2A_{n:\upsilon} = a_{n},B_{n:r} = B_{n:\upsilon} = 0.$  Moreover, if $m>0,$
 $$P({\cal R}^{(n)}_{\nu_{n}}(m,k)\leq r){\renewcommand{\arraystretch}{0,1} \begin{array}{c}
\stackrel{w}{\longrightarrow}\\
\scriptstyle{n}\end{array}}1-e^{-r^{-1}}$$ and   $$P({\cal V}^{(n)}_{\nu_{n}}(m,k)\leq v){\renewcommand{\arraystretch}{0,1} \begin{array}{c}
\stackrel{w}{\longrightarrow}\\
\scriptstyle{n}\end{array}}e^{-(-v)^{-1}},$$
respectively, with $A_{n:r} = 2A_{n:\upsilon} = c_{n},B_{n:r} = B_{n:\upsilon} = 0.$
Finally, when $-1<m<0,$ $$P(\mathcal{R}^{(n)}_{\nu_{n}}(m,k)\leq r){\renewcommand{\arraystretch}{0,1} \begin{array}{c}
\stackrel{w}{\longrightarrow}\\
\scriptstyle{n}\end{array}}\int_{0}^{\infty}\left[1-\Gamma_{\ell}({z}{r^{-(m+1)}})\right]e^{-z}dz$$ and the df of ${\cal V}^{(n)}_{\nu_{n}}(m,k)$ converge weakly to the same limit, with $A_{n:r} = 2A_{n:\upsilon} = a_{n},B_{n:r} = B_{n:\upsilon} = 0.$ \\
{\bf Example 3.2 (Pareto distribution).} It can be shown that, for the Pareto distribution $F(x) =
(1-x^{-\sigma})I_{[1,\infty)}(x), \sigma > 0.$ Therefore, in view
of Theorems 1.1 and 2.1, Part 1, in Barakat et al. (2015b), since $~~\frac{c_{n}}{a_{n}}{\renewcommand{\arraystretch}{0,1} \begin{array}{c}
\stackrel{}{\longrightarrow}\\
\scriptstyle{n}\end{array}}0,\forall ~m.$
Then$$P({\cal R}^{(n)}_{\nu_{n}}(m,k)\leq r){\renewcommand{\arraystretch}{0,1} \begin{array}{c}
\stackrel{w}{\longrightarrow}\\
\scriptstyle{n}\end{array}}\int_{0}^{\infty}\left[1-\Gamma_{\ell}(zr^{-\sigma(m+1)})\right]e^{-z}dz$$and the df of ${\cal V}^{(n)}_{\nu_{n}}(m,k)$ converge weakly to the same limit, where $A_{n:r} = 2A_{n:\upsilon} = a_{n},B_{n:r} = B_{n:\upsilon} = 0.$\\
{\bf Example 3.3 (uniform distribution).} For the uniform $(-\theta,\theta)$ distribution, by using Theorem 2.1, Part 1, in Barakat et al. (2015b),
since $\frac{a_{n}}{c_{n}}{\renewcommand{\arraystretch}{0,1} \begin{array}{c}
\stackrel{}{\longrightarrow}\\
\scriptstyle{n}\end{array}}1,$ if $m=0.$ Then
$$P(R^{(n)}_{\nu_{n}}(0,k)\leq r){\renewcommand{\arraystretch}{0,1} \begin{array}{c}
\stackrel{w}{\longrightarrow}\\
\scriptstyle{n}\end{array}}\int_{0}^{\infty}\int_{-\infty}^{0} \left[1-\Gamma_{k}(z(y-r)^{}) \right] ze^{z(y-1)}dydz$$
and $$P(V^{(n)}_{\nu_{n}}(0,k)\leq v){\renewcommand{\arraystretch}{0,1} \begin{array}{c}
\stackrel{w}{\longrightarrow}\\
\scriptstyle{n}\end{array}}\int_{0}^{\infty}\int_{-\infty}^{\infty} \left[1-\Gamma_{k}(z(y-v)^{}) \right] ze^{-z(y+1)}dydz,$$ respectively, with  $A_{n:r} = 2A_{n:\upsilon} = a_{n},B_{n:r} = B_{n:\upsilon} = 0.$\\
{\bf Example 3.4 (Beta$(\alpha, \beta)$ distribution).} For the beta distribution $F(x; \alpha,\beta), 0\leq x \leq1,\alpha,\beta>0.$
Therefore, in view of Theorem 2.2, Part 5, in Barakat et al. (2015b), if $\alpha= (m+1)\beta.$ Then
$$ P({\cal R}^{(n)}_{\nu_{n}}(m,k)\leq r){\renewcommand{\arraystretch}{0,1} \begin{array}{c}
\stackrel{w}{\longrightarrow}\\
\scriptstyle{n}\end{array}}\int_{0}^{\infty}\int_{-\infty}^{0}\left[1-\Gamma_{\ell}(z(y-r)^{\alpha})\right]z \eta \alpha (-\eta y)^{\alpha-1}e^{-z(1+(-\eta y)^{\alpha})}dydz$$ and
 $$P({\cal V}^{(n)}_{\nu_{n}}(m,k)\leq v){\renewcommand{\arraystretch}{0,1} \begin{array}{c}
\stackrel{w}{\longrightarrow}\\
\scriptstyle{n}\end{array}}\int_{0}^{\infty}\int_{-\infty}^{\infty}\left[1-\Gamma_{\ell}(z(y-v)^{\alpha})\right]z \eta \alpha (\eta y)^{\alpha-1}e^{-z(1+(\eta y)^{\alpha})}dydz,$$
respectively, where $\eta= (\frac{\beta}{c})^{\frac{1}{\beta}}(\frac{m+1}{c\alpha})^{\frac{1}{\alpha}},c=\frac{\Gamma(\alpha+\beta)}{\Gamma(\alpha)\Gamma(\beta)}.$
 Clearly, the same result holds for the power distribution $F(x;\alpha, 1).$\\
{\bf Example 3.5 (standard normal, logistic, Laplace, and log-normal distributions).}
 After some algebra, we get,$$\eta^{-1}=\lim_{n \to \infty}\frac{c_{n}}{a_{n}}=\!\!\left\{
    \begin{array}{ll} \frac{1}{\sqrt{m+1}}, ~~~\mbox{for the normal distribution,}
\\1,~~~~~~~ \mbox{for the logistic and Laplace distribution,}
\\0, ~~~~~~~\mbox{for the log-normal distribution.}
\end{array}\right.$$
Moreover, for the standard normal, logistic, and Laplace distributions, we get
$$P({\cal R}^{(n)}_{\nu_{n}}(m,k)\leq r)
{\renewcommand{\arraystretch}{0,1} \begin{array}{c} \stackrel{w}{\longrightarrow}\\ \scriptstyle{n}\end{array}}\int_{0}^{\infty}\int_{-\infty}^{\infty}\left[1-\Gamma_{\ell}(ze^{(y-r)(m+1)})\right]z \eta e^{-\eta y}e^{-z(1+e^{-\eta y})}dydz$$
and $$P({\cal V}^{(n)}_{\nu_{n}}(m,k)\leq v){\renewcommand{\arraystretch}{0,1} \begin{array}{c}
\stackrel{w}{\longrightarrow}\\
\scriptstyle{n}\end{array}}\int_{0}^{\infty}\int_{-\infty}^{\infty}\left[1-\Gamma_{\ell}(ze^{(y-v)(m+1)})\right]z \eta e^{\eta y}e^{-z(1+e^{\eta y})}dydz,$$
for the standard normal distribution  $(m=0,k=1).$ Then
$$~P({\cal R}_{\nu_{n}}^{}(0,1)\leq a_{n}r+2b_{n}){\renewcommand{\arraystretch}{0,1} \begin{array}{c} \stackrel{w}{\longrightarrow}\\ \scriptstyle{n}\end{array}} \int_{0}^{\infty} \frac{y^{2}~dy}{(y^{2}+y+e^{-r})^2} =\left\{
\begin{array}{ccc}
f_{1}(r), ~~~~~~~~~~~& \mbox{~$r<\ln{4}$}, &\\
\frac{2}{3}, ~~~~~~~~~~~~& \mbox{ $r$=$\ln{4}$}, &  \\
f_{2}(r), ~~~~~~~~~~~~~~& \mbox{ $r>\ln{4}$}, &
\end{array}
\right.$$
where
$$~f_{1}(r)=(4e^{-r}-1)^{-1}\left[\frac{4e^{-r}}{\sqrt{4e^{-r}-1}}~cot^{-1}~(\frac{1}{\sqrt{
4e^{-r}-1}})-1\right] ,$$ $$~f_{2}(r)=(1-4e^{-r})^{-1}\left[1-\frac{2e^{-r}}{\sqrt{1-4e^{-r}}}\ln{\frac{1+\sqrt{
1-4e^{-r}}}{1-\sqrt{1-4e^{-r}}}}\right],$$ $$~a_{{n}}=\frac{1}{\sqrt{2\ln{n}}},~b_{n}=\sqrt{2\ln{n}}~~- ~~ \frac{\ln{\ln{n}}+\ln{4\pi}
}{2\sqrt{2\ln{n}}}~$$ and
$$P({\cal V}_{\nu_{n}}(0,1)\leq a_{n}v)
 {\renewcommand{\arraystretch}{0,1} \begin{array}{c} \stackrel{w}{\longrightarrow}\\ \scriptstyle{n}\end{array}}
 1-\int_{0}^{\infty} \frac{dy}{(y(e^{2v}+1)+1)^{2}}=(1+e^{-2v})^{-1}.$$
Finally, for the log-normal distribution,\\ $P({\cal R}^{(n)}_{\nu_{n}}(m,k)\leq r){\renewcommand{\arraystretch}{0,1} \begin{array}{c}
\stackrel{w}{\longrightarrow}\\
\scriptstyle{n}\end{array}}\int_{0}^{\infty} \left[1-\Gamma_{\ell}(ze^{-r(m+1)}) \right]e^{-z}dz$ and the df of ${\cal V}^{(n)}_{\nu_{n}}(m,k)$ converge weakly to the same limit. \\
{\bf Example 3.6 (exponential $(\sigma)$ and Rayleigh $(\sigma)$ distributions).} In view of Theorem 2.2, Part 4 in Barakat et al. (2015b). Then\\
$P({\cal R}^{(n)}_{\nu_{n}}(m,k)\leq r) {\renewcommand{\arraystretch}{0,1} \begin{array}{c}
\stackrel{w}{\longrightarrow}\\
\scriptstyle{n}\end{array}}\int_{0}^{\infty}\left[1-\Gamma_{\ell}(ze^{-r(m+1)})\right]
e^{-z}dz$ and the df of ${\cal V}^{(n)}_{\nu_{n}}(m,k)$ converge weakly to the same limit,
 for exponential $(\sigma)$ and Rayleigh $(\sigma)$ distributions.

\noindent {\bf Acknowledgements}

Elsawah's work was partially supported by the UIC GRANT R201409 and the Zhuhai Premier
Discipline Grant and Qin's work was partially supported by the National Natural Science Foundation of China (Nos. 11271147, 11471135, 11471136).

\end{document}